
\input amstex
\documentstyle{amsppt}
\refstyle{A}
\NoRunningHeads

\topmatter
\title On a theorem of Harer \endtitle
\author Robert  Treger \endauthor
\address AG{\&}N, Princeton, New Jersey  \endaddress
\email tregrob{\@}aol.com \endemail
\endtopmatter

\document
\head Introduction \endhead

Let  ${\bold M}_g$ be the coarse moduli space of smooth projective complex
curves of genus $g$.  Let  $\text{\rm Pic} ({\bold M}_g)$
denote the Picard group of line bundles on ${\bold M}_g$.  According to
Mumford,  $\text{\rm Pic} ({\bold M}_g) \otimes \Bbb Q \simeq H^2 ({\bold
M}_g , \Bbb Q )$  for $g \geq 2$ \cite{9}.

The purpose of this note is to give another proof of the following
fundamental theorem of Harer \cite{4}.
\proclaim {Theorem 1} $\text{\rm rk}\  \text{\rm Pic} ({\bold M}_g)
\otimes \Bbb Q =1$  for $g \geq 3.$
\endproclaim
We will deduce Theorem 1 from the Harer stability theorem \cite{5, 6, 8}.
\proclaim {Theorem 2} $H^2 ({\bold M}_g , \Bbb Q) \simeq H^2 ({\bold M}_{g+1}
, \Bbb Q)$  for $g \geq 6.$
\endproclaim
Note that Theorem 2 is similar to the corresponding stability results for
arithmetic groups.

\head Proof of Theorem 1 \endhead

To deduce Theorem 1 from Theorem 2, it suffices to show that
$\text{\rm rk} \ \text{\rm Pic} ({\bold M}_g)
\otimes \Bbb Q =1$  for $g =3, 4, 5, 6$.
Let $V_{n,d}\subset {\Bbb P}^N \  (N = n(n+3)/2)$
be the variety of irreducible complex plane curves of degree  $n$ with $d$
nodes and no other singularities.  Let $V(n,g)\subset {\Bbb P}^N \  (N =
n(n+3)/2)$
be the variety of irreducible complex plane curves of degree  $n$ and genus
$g =(n-1)(n-2)/2 - d$.  Clearly $V_{n,d}\subset V(n,g).$
Let
$$\Sigma_{n,d} \subset {\Bbb P}^N \times {\text
{Sym}}^d({\Bbb P}^2)
$$
be the {\it closure\/} of the locus of pairs  $(E,\, \sum_{i=1}^dR_i),$
where  $E$  is an irreducible
nodal curve and  $R_{1}$, \dots , $R_{d}$  are its nodes,  and $\pi_{n,d}$
the projection to the symmetric product.
We fix a general point  $P$  and a general line $L$ in ${\Bbb P}^2 ,$ and
consider the following two elements in
$\text{\rm Pic} (V_{n,d})$: $(CP) =$ the divisor class of curves containing
the point $P$, and $(NL)=$ the divisor class of curves with a node located
somewhere on  $L$.

If $g\geq 3(g+2-n)$ then the natural map $V(n,g)  \rightarrow {\bold M}_g$ is
a {\it surjective \/} morphism \cite{3, p. ~358}.  Precisely, for any smooth
curve $C$ of genus $g$, there is a line bundle $\Cal L$ on $C$ of degree $n$
such that $h^0 (C, \Cal L ) \geq 3$.  Furthermore, Harris observed that
 $\text{\rm rk}\  \text{\rm Pic} ({\bold M}_g)
\otimes \Bbb Q =1$  provided one can show that $\text{\rm Pic} (V_{n,d})$ is
a torsion group, where $d =(n-1)(n-2)/2 - g$ \cite{7}.  Thus Theorem 1
follows from the
\proclaim{Lemma} Let $n$ and $d$ be a pair of positive integers such that one
of the following conditions holds: i) $n \gg d$ or ii) $(n,d) = (5,3), (6,6),
(6,5), \text{or} \ (6,4)$.    Then $\text{\rm Pic} (V_{n,d})$ is a torsion
group generated by $(CP) $ and $(NL)$.
\endproclaim

\demo{Proof of Lemma}
We will prove the lemma for $(n,d) = (6,6)$, hence $g=4$.  The remaining
cases are similar only easier. Given {\it any\/} $6$ distinct points $q_1,$
\dots , $q_{6}$ in ${\Bbb P}^{2}$, no four on a line, then a general curve
$C$ of the linear system $\Cal L$ of curves of degree $6$ with assigned
singularities at $q_1,$ \dots , $q_{6}$ has  $6$  nodes and no other
singularities.  We will calculate $\dim {\Cal L}$ using Castelnuovo's method
(see \cite {1, Book IV, Chap.~I, Sect.~3}).

The characteristic series determined on $C$ by $\Cal L$ is cut out by curves
of $\Cal L$.  Clearly one can find two curves in $\Cal L$, say $C_1$ and
$C_2$, such that  $(C_1 \cdot C_2)_{q_i} = 4$ for $i= 1,$ \dots , $6.$ It
follows that the degree of the characteristic series equals ${(\deg C)}^2 - 6
\cdot 4 > 6 = 2g - 2.$  Hence this series is non-special, the superabundance
of $\Cal L$ vanishes, and $\dim {\Cal L} = 6(6+3)/3 - 6 \cdot 3.$

Let $S$ be a closed subset of ${\text
{Sym}}^6({\Bbb P}^2)$ which is a union of the singular locus and the cycles
with at least $4$ points on a line.  Clearly codim $\! S =2$ and $\text{\rm
Pic} ({\text  {Sym}}^6({\Bbb P}^2 ) \setminus S) =\text{\rm Pic} ({\text
 {Sym}}^6({\Bbb P}^2 ) ) = {\Bbb Z}.$ We set
$$
\Sigma = {\pi}_{6,6}^{-1} ({\text  {Sym}}^6({\Bbb P}^2 ) \setminus S) \cap
\Sigma_{6,6}.
$$
 By the above discussion, all the fibers of ${\pi}_{6,6} | \Sigma$ have the
same dimension and ${\pi}_{6,6} | \Sigma$: \ $\Sigma \rightarrow {\text
 {Sym}}^6({\Bbb P}^2 ) \backslash S$ is a projective bundle.  Thus $\text{\rm
Pic} (\Sigma) = {\Bbb Z}^2$ and its generators correspond to the divisor
classes $(CP)$ and $(NL)$ coming from the fiber and the base of the
fibration, respectively.

We consider two divisors on $\Sigma$: (i) the closure of the locus of
irreducible curves with  $d - 1$  assigned nodes and one assigned cusp,  and
 (ii) the closure of the locus of nodal curves with  $d$  assigned nodes and
one unassigned node.  The corresponding elements of  $\text{\rm Pic} (\Sigma)
\otimes \Bbb Q$  are linearly independent \cite{2, Sect.~2}. There is a
natural open immersion  $V_{6,6} \subset \Sigma .$  Hence we get an exact
sequence
$$
{\Bbb Z}^r  \rightarrow  \text{\rm Pic} (\Sigma)  \rightarrow  \text{\rm Pic}
(V_{6,6})  \rightarrow 0
$$
where  ${\Bbb Z}^r$  is generated by the irreducible components of  $\Sigma
\setminus  V_{6,6}$  of codimension one  in  $\Sigma .$  Therefore
 $\text{\rm Pic} (V_{6,6})$  is torsion.
This proves the lemma and the theorem.
\enddemo

\enddocument